\documentclass{amsart}

\usepackage{amsfonts}
\usepackage{amsthm}
\usepackage{amsmath}
\usepackage{amsfonts}
\usepackage{latexsym}
\usepackage{amssymb}
\usepackage[latin1]{inputenc}
\newcommand{\lin}{\text{span}}

\def\a{\mathbf a}

\def\K{{\mathcal K}}

\newcommand{\U}{\mathcal{U}}

\def\rai{^{1/2}}

\def\barr{\begin{array}}
\def\earr{\end{array}}

\def\a{\mathbf{a}}
\def\b{\mathbf{b}}
\def\c{\mathbf{c}}
\def\v{\mathbf{v}}
\def\e{\mathbf{e}}
\newcommand{\CC}{\mathbb{C}}
\newcommand{\RR}{\mathbb{R}}
\newcommand{\NN}{\mathbb{N}}

\newcommand{\FF}{\mathbb{F}}
\def\NN{\mathbb{N}}

\def\C{\mathbb{C}}

\def\ra{\rightarrow}
\def\bma{\left[\begin{array}}
\def\ema{\end{array}\right]}
\def\ben{\begin{enumerate}}
\def\een{\end{enumerate}}

\newtheorem{fed}{Definition}[section]
\newtheorem{teo}[fed]{Theorem}
\newtheorem*{teo*}{Theorem}
\newtheorem{lem}[fed]{Lemma}
\newtheorem{cor}[fed]{Corollary}
\newtheorem{pro}[fed]{Proposition}
\theoremstyle{definition}
\newtheorem{rem}[fed]{Remark}

\newtheorem{exa}[fed]{Example}

\def\eps{\varepsilon}

\def\la{\lambda}

\def\cA{\mathcal{A}}
\def\cB{\mathcal{B}}

\def\cD{\mathcal{D}}
\def\cE{\mathcal{E}}

\def\cF{\mathcal{F}}
\def\cK{\mathcal{K}}

\def\cM{\mathcal{M}}

\def\cP{\mathcal{P}}

\def\cR{\mathcal{R}}

\def\cT{\mathcal{T}}

\def\cV{\mathcal{V}}

\def\cG{\mathcal{G}}

 \DeclareMathOperator{\tr}{tr}

\DeclareMathOperator{\diag}{diag}

\DeclareMathOperator{\sspan}{span}

\newcommand{\ds}{\displaystyle}

\newcommand{\pint}[1]{\displaystyle \left \langle #1 \right\rangle}

\newcommand{\hil}{\mathcal{H}}

\newcommand{\norm}[1]{\left\| #1 \right\|}

\begin{document}

\title{Minimization of convex functionals over frame operators}
\author{ Pedro Massey}
\address{Dpto. de Matemática, Univ. Nac. de La Plata and IAM-CONICET}
\email{massey@mate.unlp.edu.ar}

 \author{ Mariano Ruiz}
 \address{Dpto. de Matemática, Univ. Nac. de La Plata and IAM-CONICET}
\email{mruiz@mate.unlp.edu.ar}
\thanks{Partially supported by CONICET (PIP 5272) and
UNLP (11 X350)}

\begin{abstract}  We present results about minimization of convex functionals defined over   a finite set of vectors in a finite dimensional Hilbert space, that extend several known results for the Benedetto-Fickus frame potential. Our approach depends on majorization techniques. We also consider some perturbation problems, where a positive perturbation of the frame operator of a set of vectors is realized as the frame operator of a set of vectors which is close to the original one.
\end{abstract}
\maketitle

\section{Introduction}
Let $\hil$ be a Hilbert space. A set of vectors $\cF=\{\phi_i\}_{i\in I}$ in $\hil$  is a {\sl frame} if there exist a pair of constants $a,b>0$ such that, for every $x\in \hil$,
\begin{equation}\label{desig frame intro}  a\ \|x\|^2\leq \sum_{i\in I}| \pint{x \, ,\,  \phi_i}|^2\leq b\ \|x\|^2.
\end{equation} 
The optimal constants  $a,b$ in \eqref{desig frame intro} are called the {\sl frame bounds}. We say that the frame is {\sl tight} if $a=b$. In general, if the inequality on the right hand side of \eqref{desig frame intro} holds for $x\in \hil$ we say that $\cF$ is a Bessel sequence. Given a Bessel sequence $\cF$ we consider its {\sl synthesis operator} $T^\cF:l_2(I)\ra \hil$ defined as $T^\cF(e_i)=\phi_i$, where $\{e_i\}_{i\in I}$ is the canonical orthonormal basis of $l_2(I)$. We also consider its {\sl frame operator} given by $S^\cF=T^\cF(T^\cF)^*$ and its {\sl Grammian}, defined by $G^\cF=(T^\cF)^*T^\cF$.

Frames where introduced by Duffin and Schaeffer \cite{DS} in their work on nonharmonic Fourier series. These were later rediscovered by Daubechies, Grossmann and Meyer in the fundamental paper \cite{[DGM]}. In recent years the study of frames has  increased considerably  due to the wide range of applications in which frames play an important role. In this note we shall focus on finite frames i.e. $\hil= \mathbb{F}^d$ where $\mathbb{F}=\CC$ or $\RR$ and $I$ is a finite set. Note that in this setting, a frame is just a set of generators for $\hil$.

In \cite{BF} Benedetto and Fickus introduced the notions frame force (FF) and frame potential (FP) for a finite frame. More explicitly they defined, for $\cF=\{\phi_i\}_{i=1}^m\subseteq \hil$ a finite sequence of vectors \begin{equation}\label{intro fp}\mathrm{FP}(\cF)=\sum_{i,j=1}^m|\langle \phi_i,\phi_j\rangle|^2=\tr((S^\cF)^2)\end{equation}
It is shown in \cite{BF} that the finite unit norm tight frames are the minimizers of the frame potential among all unit norm frames with a fixed number of vectors. If we now impose restrictions on the lengths of the vectors, the structure of minimizers changes since tight frames with a prescribed set of norms may not exist. The complete characterization of global and local minimizers for the frame potential was done in \cite{casazza2}.

The equality $\mathrm{FP}(\cF)=\tr((S^\cF)^2)$ suggests that, more generally,  we can consider functionals of the form $P_f(\cF)=\tr(f(S^\cF))$, where $f$ is a non-negative, non-decreasing and convex function defined on $[0,\infty)$. In this context, the problem of describing the geometrical structure of minimizers of these convex functionals arises; surprisingly, this structure  does not depend on $f$. In order to state the following results we introduce the sets 
$\cA(c)=\{\{\phi_i\}_{i=1}^m\subset \CC^d,\,\sum_{i=1}^m\|\phi_i\|^2=c\}$ and
$\cB(\a)=\{\{\phi_i\}_{i=1}^m\subset\CC^d,\, \|\phi_i\|^2=a_i \text{ for every } i\}$, where $\a=(a_i)_{i=1}^m$ is a non-increasing finite sequence of positive real numbers.

\begin{teo*}[A]
Let  $f:\RR_{\geqslant 0}\rightarrow \RR_{\geqslant 0}$ be a convex function and $P_f$ the functional associated to $f$. Let $c>0$ and $\a=(a_i)_{i=1}^m$ be a non-decreasing finite sequence of positive real numbers. Then, 
\begin{itemize}
\item[a)] If $\cF\in \cA(c)$ is a tight frame then it is a global minimizer of $P_f$ in $\cA(c)$. If we assume further that $f$ is strictly convex then every global minimizer in $\cA(c)$ is  tight. 
\item[b)] 
If $\cF \in \cB(\a)$ is of the form
\begin{equation}\label{intro minimos globales}
 \{ \sqrt{a_i} \,
e_i\}_{i=1}^r\cup \{\phi_i\}_{i=r+1}^m
\end{equation}
where $\{e_i\}_{i=1}^d$
is an o.n.b. for $\CC^d$, $r$ is the
\emph{d-irregularity} of $\a$ (see definition \ref{defi irre} below) and
  $\{\phi_i\}_{i=r+1}^m$ is a
tight frame for $\text{span }\{e_i\}_{i=r+1}^d$ then, it is a global minimizer of $P_f$.
If we assume further that $f$ is strictly convex then every global minimizer in $\cB(\a)$ is as in \eqref{intro minimos globales} for some o.n.b. $\{e_i\}_{i=1}^d$.
\end{itemize}
\end{teo*}

It is also interesting to study the structure of the \emph{local} minimizers of $P_f$ in the previous sets $\cA(c)$ and $\cB(\a)$. A natural metric in this context is the vector-vector distance $d(\cF,\cG)=\max_{1\leq i\leq m}\|\phi_i-\psi_i\|$ for sequences $\cF=\{\phi_i\}_{i=1}^m$, $\cG=\{\psi_i\}_{i=1}^m$. But this characterization problem turns out to be quite difficult for the local minimizers of $P_f$ in $\cB(\a)$. Hence, we alternatively consider the description of the structure of local minimizers of $P_f$ in $\cR(\a)=\{S^\cF, \ \cF\in \cB(\a) \}$ endowed with the norm topology. Notice that this last point of view is weaker. Indeed, $\|S^\cF-S^\cG\|\leq 2\sqrt{m}\,\max(\| T^\cF\|, \|T^G\|) \,d(\cF,\cG)$, where $T^\cF$ and $S^\cF$ denote the synthesis and frame operator of $\cF$ (see the beginning of section \ref{prelims frames}), while there are pairs of different sequences that share the frame operator.
\begin{teo*}[B] 
Let  $f:\RR_{\geqslant 0}\rightarrow \RR_{\geqslant 0}$ be a non decreasing strictly convex function and $P_f$ the functional associated to $f$. Let $c>0$ and $\a=(a_i)_{i=1}^m$ be a non-decreasing finite sequence of positive real numbers. Then, 
\begin{itemize}
\item[a)] Every local minimizer of $P_f$ in $\cA(c)$ with respect to $d(\cdot,\cdot)$ is a tight frame and hence, a global minimizer.
\item[b)] Every local minimizer of $P_f$ in $\cR(\a)$ with respect to the operator norm is of the form \eqref{intro minimos globales} for some o.n.b. $\{e_i\}_{i=1}^d$ of $\CC^d$ and hence, a global minimizer.
\end{itemize}
\end{teo*}

The previous results show that the structure of the local minimizers of $P_f$ (when $P_f$ is considered as a function of the frame operators) does not depend on the strictly convex function chosen.  Unfortunately, we get only partial results related with the local minimizers of $P_f$ in $\cB(\a)$ with respect to the vector-vector distance, for a general convex function $f$. 

Our approach depends on solving some
perturbation problems concerning the frame operator for a generic case of frame.

More explicitly, 
if $\cF$ is a frame in $\cB(\a)$ which can not be partitioned in two mutually orthogonal sets of vectors (i.e. its Grammian  is not block-diagonal) and $S_i$ is a sequence in $\cM_d(\CC)^+$ which converges to $S^\cF$, then for every $\eps>0$ there exists $i_0$ such that, for $i\geq i_0$ there is a frame $\cG \in \cB(\a)$ such that $S^\cG=S_i$ and 
$d(\cF,\cG)\leq \eps$. Our approach to this problem depends on differential geometric tools that we describe in an appendix at the end of the paper.  
In the particular case of the Benedetto-Fickus frame potential, we recover a theorem by Casazza et al. \cite{casazza2} describing its local minimizers.

The paper is organized as follows:  Section \ref{Preliminares} contains preliminary facts together with some new results about majorization of vectors in $\RR^d$ that we shall need in the sequel; Propositions \ref{maximos y minimos} and \ref{hay min para la mayo} give a characterization of minimal points of certain sets of vectors with respect to majorization.
Section \ref{prelims frames} is devoted to  the basic facts about frames in $\CC^d$ together with some previous results from \cite{JDM} about some design problems for frames. In Section \ref{Convex funct}, some properties of the convex functions $P_f$ defined on frame operators are given. In this section we consider  the sets of frame operators $\cR(c)$ and $\cT(\a)$, consisting of frame operators of elements in $\cA(c)$ and $\cB(\a)$ respectively.  Theorems  \ref{sin normas} and \ref{T:minimos globales} deal with the characterization of global and local minimizers for every $P_f$ (for a non decreasing strictly convex function $f:\RR_{\geqslant 0}\rightarrow \RR_{\geqslant 0}$) on $\cR(c)$ and $\cT(\a)$. At the end of this section, some examples and applications are given.
Finally, in Section \ref{frame op to frames}  we focus on the structure of minimizers of the functions $P_f$ when they are defined on frames instead of frame operators. This leads to some geometrical problems which are developed in the Appendix.

\medskip

\noindent {\bf Acknowledgments.} We would like to thank Demetrio Stojanoff and Jorge Antezana for several useful comments related with the content of this note that helped us improve its exposition.

\section{Preliminares}\label{Preliminares}

In this section we present some basic aspects of majorization
theory together with some new results that we shall need in what follows.
 For a more detailed treatment of majorization see
\cite{horn2}.
Given $\b=(b_1,\ldots,b_d)\in \RR^d$, denote by  $\b^ \downarrow
\in \RR^d $ the vector obtained by rearranging the coordinates of
$\b$ in non increasing order. If $\b,\, \c\in \RR^d$ then we say
that $\b$ is \emph{majorized} by $\c$, and write $\b\prec  \c$,
if $$ \sum_{i=1}^k b^\downarrow_i\leq \sum_{i=1}^k c^\downarrow_i
\ \ \ \ k = 1,\dots,d-1\  \text{ and } \;\sum_{i=1}^ d b^\downarrow _i=\sum_{i=1}^ d c^\downarrow _i. $$

Majorization is a preorder relation in $\RR^d$ that occurs
naturally in matrix analysis and plays an important role in convex
optimization problems.

\begin{pro}\label{maximos y minimos}
Let $c>0$ and consider the set \begin{equation}\label{defi de k y
c } \K(c)=\{\b\in (\RR_{\geqslant 0})^d:\ \sum_{i=1}^db_i=c\}
\end{equation} Then the vector
$\v=(\frac{c}{d},\ldots,\frac{c}{d})$
 satisfies $\b\succ \v$ for every  $\b\in \K(c)$.
Moreover, if $\b\in \K(c)$ is such that $\b^\downarrow\neq \v$,
 then  for every $0<\eps $ sufficiently small, there exists $\b_\eps \in \K(c)$ such that $ \b^\downarrow\neq\b_\eps^\downarrow$, $\b\succ \b_\eps$ and $ \| \b^\downarrow - \b_\eps^{\downarrow} \|\leq  \eps$.

\end{pro}
\begin{proof}
The first part is a well known fact about majorization, and it is easy to check.  For the proof of the moreover part, suppose that $\b \in \K(c)$ is such that $\b\neq \v$, then there exists  a index $j$, $1 \leq j \leq d$ such that $b_j^\downarrow > b_{j+1}^\downarrow$ were we denote by $b_i^\downarrow$ the entries of $\b^\downarrow$.

Let $0<\eps$ such that $b_j^\downarrow -\sqrt{\frac{\eps}{2}}\ \geq b_{j+1}^\downarrow +\sqrt{\frac{\eps}{2}}$  and denote by $\b_\eps$ the vector  $\b_\eps =\b^\downarrow - \sqrt{\frac{\eps}{2}}\,\e_j +\sqrt{\frac{\eps}{2}}\,\e_{j+1}$ were $\{\e_i\}_{i=1}^d $ is the  canonical basis in $\RR ^d$.
Clearly $\b_\eps \in \K(c)$, $\b \succ \b_\eps$, and  by construction of $\b_\eps$,
$\|\b^\downarrow - \b_{\eps}^{\downarrow}\|^2= \eps$.
\end{proof}

Following \cite{casazza2} we consider the $d$-irregularity of a sequence as follows
\begin{fed}\label{defi irre} \rm Let $\a=(a_i)_{i=1}^m$ be a non increasing sequence of positive numbers  and
$d \in \NN$ with $d\leq m$. The $d$-irregularity of $\a$, denoted  $r_d(\a) \in \NN$,  is defined as
$$r_d(\a) = \max \Big\{ 1\leq j\leq d-1  :
(d-j)a_j>\sum_{i=j+1}^m a_i\, \Big\},$$
if the set on the right is non empty, and $r_d(\a)=0$ otherwise.
\end{fed}
Notice that in particular, with the notations of Definition \ref{defi irre}, we have: 
\begin{enumerate}
\item $(d-j)a_j\leq \sum_{i=j+1}^m a_i$,  for $r_d(\a)< j\leq d$  whenever $r_d(\a)>0$,
\item  $(d-j)a_j>\sum_{i=j+1}^m a_i\, $,   for every $1\leq j\leq r_d(\a)$.
\end{enumerate}

\begin{pro}\label{hay min para la mayo}
Let $0<d\leq m$ and let $\a=(a_i)_{i=1}^m$ be a non increasing  sequence of positive numbers
with $d$-irregularity $r=r_d(\a)$. Consider the set

$$\cP(\a)=\{  \b\in (\RR_{\geqslant 0})^d : \  \sum_{i=1}^k b_i^\downarrow \geq \sum_{i=1}^k a_i \;  \text { for } \;  1\leq k\leq d  \; \text { and }\;  \sum_{i=1}^d b_i = \sum_{i=1}^m a_i\}.$$
Let
$\v=(a_1,\ldots,a_r, \overbrace{c, \ldots, c}^{d-r\;  times})$, where $c=(d-r)^{-1}\sum_{j=r+1}^m a_j$.
Then $\v$ belongs to $\cP(\a)$ and, for every $\b \in \cP(\a)$, $\b \succ \v$.
Moreover, if $\b \in \cP(\a)$ and $\b^\downarrow\neq \v$, then  for every $0<\eps $ sufficiently small,  there exists $\b_{\eps}$ in $\cP(\a)$ such that $\b_{\eps}^\downarrow\neq \b^\downarrow$, $\b \succ \b_{\eps}$ and $\|\b^\downarrow - \b_\eps^{\downarrow} \|\leq \eps$.

\end{pro}
\begin{proof}
By the comments after Definition \ref{defi irre}, $\v= \v^\downarrow$. First, we show
that $\v\in \cP(\a)$. Note that $\sum_{j=1}^k a_j=\sum_{j=1}^k v_j$ for
$1\leq j\leq r$. On the other hand,
$$a_{r+1}(d-r)-a_{r+1}=a_{r+1}(d-(r+1))\leq
\sum_{j=r+2}^m a_j \ \Rightarrow \
a_{r+1}\leq(d-r)^{-1}\sum_{j=r+1}^m a_j.$$ Therefore
$c\geq a_{r+1}\geq a_{j}$ for every $r+1\leq j\leq m$. Then,
for every $r+1\leq k\leq d$ we have $$\sum_{j=1}^k
v_j=\sum_{j=1}^r a_j+\sum_{j=r+1}^k c_j\geq \sum_{j=1}^k
a_j.$$ Since $\sum_{j=1}^d v_j=\sum_{j=1}^m a_j$ it follows
that $\v\in \cP(\a)$.
 Let $\b=(b_i)_{i=1}^d \in \cP(\a)$ and, without loss of generality, assume that $\b=\b^\downarrow$. Then, it is clear that $\sum_{j=1}^k v_j\leq
\sum_{j=1}^k b_j$ for every $1\leq k\leq r$. Let
$\alpha=\sum_{j=1}^r b_j-\sum_{j=1}^r a_j \geq 0$. Therefore
\begin{equation} \left(\sum_{j=1}^r b_j-\sum_{j=1}^r a_j
\right)+ \sum_{j=r+1}^d b_j= \sum_{j=r+1}^m a_j \
\Rightarrow \sum_{j=r+1}^d( b_j+(d-r)^{-1}\alpha)=
\sum_{j=r+1}^m a_j
\end{equation} which implies, by Proposition \ref{maximos y minimos},
that
$(c)_{i=r+1}^d\prec((d-r)^{-1}\alpha+ b_i)_{i=r+1}^d\in\RR^{d-r}$.
Then, for every $r+1\leq k\leq d$ we have
\begin{eqnarray*}
\sum_{j=1}^k b_j&=&\sum_{j=1}^r b_j-\sum_{j=r+1}^k
(d-r)^{-1}\alpha +\sum_{j=r+1}^k(b_j+(d-r)^{-1}\alpha)\\
&\geq & \sum_{j=1}^r b_j-\alpha+\sum_{j=r+1}^k
c=\sum_{j=1}^ra_j+\sum_{j=r+1}^k c=\sum_{j=1}^k v_j.
\end{eqnarray*} 
On the other hand $$ \sum_{j=1}^d b_j=\sum_{j=1}^m a_j=\sum_{j=1}^d v_j$$
so we see that $\v\prec \b$. For the second part, let $\b \in \cP(\a)$, $\b^\downarrow \neq \v$. Again we assume that $\b=\b^\downarrow$.

\medskip

\textbf{Claim:} There exists $j$, $1\leq j \leq d-1$ such that $b_j>b_{j+1}$ and $\sum_{i=1}^j b_i > \sum_{i=1}^j a_i$.

\medskip

It is clear that for some $1\leq k \leq d-1$, $b_k>b_{k+1}$. Otherwise, $b_i=b_1$ for all $i$ which would imply that $\b=\v$ (the $d$-irregularity of $a$ would be 0).
Denote by  $b_{t_1} \geq b_{t_2} \geq \ldots \geq b_{t_m}$ all  the entries of $\b $ which satisfy $b_{t_n} > b_{t_n +1}$. 

Suppose that, for every $t_n$, $\sum_{i=1}^{t_n} b_i =\sum_{i=1}^{t_n} a_i$. Then, since by hypothesis $k b_1=\sum_{i=1}^k b_i\geq \sum_{i=1}^k a_i$ for all $k\leq t_1$, we have that $a_i=b_1=b_i$ for all $i\leq t_1$. By the same reasoning, $a_i=b_{t_1+1}=b_i$ for all $t_1+1 \leq i \leq t_2$. Finally, we get that $a_i=b_i $ for all $1\leq i \leq t_m$  moreover,  $b_k=(d-t_m)^{-1}\sum_{i=t_m+1}^m a_i$ for $t_m+1 \leq k$. The definition of the irregularity of $\a$ implies that $t_m\leq r$ (otherwise, the decreasing order of $\b$ would be violated), but if $t_m\leq r-1$, then  by the comments following Def. \ref{defi irre}, $$a_{t_m+1}>(d-(t_m+1))^{-1}\sum_{i=t_m+2}^m a_i,$$
 which in turn implies that $a_{t_m+1}> (d-t_m)^{-1}\sum_{i=t_m+1}^m a_i=b_{t_m+1}$, which contradicts $\b \in \cP(\a)$. The only possible case is $t_m=r$, but in this case, $\b=\v$, a contradiction.

\medskip

Now, given $1\leq j \leq d-1$ such that $b_j>b_{j+1}$ and $\sum_{i=1}^j \b_i > \sum_{i=1}^j a_i$, let $\eps$ such that  $b_j -{\eps}/{\sqrt{2}}  \geq b_{j+1}+ {\eps}/{\sqrt{2}}$ and  $\sum_{i=1}^j b_i - {\eps}/{\sqrt{2}}\geq   \sum_{i=1}^j a_i$. Now, denote by $\b_{\eps}$ the vector $\b-{\eps}/{\sqrt{2}}\,\e_j +{\eps}/{\sqrt{2}}\,\e_{j+1}$.
Then is easy to see that $\b_{\eps}$  satisfy the desired properties.

\end{proof}

\begin{rem}\label{ordenadito}
Note that the proof of the previous claim shows that the only vector $\b$ in $\cP(\a)$ such that: $\b^\downarrow=(a_1,a_2,\ldots,a_k,c,\ldots,c)$ is $\v$.

\end{rem}

Finally, we consider the following extension
of majorization to self-adjoint operators due to Ando \cite{An} which will be useful for the study of convex functions on frame operators: given self-adjoint matrices $B,\,C\in \cM_d(\CC)$ we
say that $B$ is majorized by $C$, and write $B\prec C$ if and
only if $\lambda(B)\prec \lambda(C)$, where $\lambda(A)\in \RR^d$
denotes the $d$-tuple of eigenvalues of a selfadjoint matrix $A\in
\cM_d(\CC)$ counted with
multiplicity and arranged in decreasing order.

\section{Preliminaries on frames}\label{prelims frames}
Let $\hil=\FF^d$, ($\FF=\CC$ or $\RR$), and let $\cF=\{\phi_i\}_{i=1}^m$  
 be a set of vectors in $\hil$, we say that $\cF$ is a {\it frame} if there exist $a,b>0$ such that for every vector $x$ in $\hil$

\begin{equation}\label{desigualdad frame}
a\|x\|^2\leq \sum_{i=1}^m |\pint{x \, , \, \phi_i}|^2 \leq b\|x\|^2
\end{equation}

the optimal bounds $a$ and $b$ are the upper and lower {\it frame bounds} for $\cF$.

We can define the following bounded linear operator $$T^\cF: \FF^m \rightarrow \hil,\ \ \
T^\cF(e_i)=\phi_i,\ 1\leq i\leq m$$
The positive semidefinite operators $$ G^\cF:=(T^\cF)^*T^\cF  \\\ \mbox{ and } \ \ \ S^\cF:=T^\cF(T^\cF)^*$$ 
are called \emph{Grammian} and  the \emph{frame operator} respectively,  of the  sequence $\cF=\{\phi_i\}_{i=1}^m$.
Throughout this note we shall consider the matrices of those operators with respect to the canonical bases of $\FF^m$ and $\FF^d$, maintaining the notation . Thus, $S^\cF\in \mathcal{M}_d(\FF)^+ $ and $G^\cF \in\mathcal{M}_m(\FF)^+$.

In particular, it can be seen that the upper and lower frame bound for $\cF$ are  the greatest and smallest positive eigenvalues of $S^\cF$, denoted by $\la_1$ and $\la_d$ respectively.

\begin{pro}[\cite{JDM}]\label{el s y el g}
Let $\cF=\{\phi_i\}_{i=1}^m\subseteq \hil$ and let $G$ and $S$
be the Grammian and frame operators of $\cF$. Then, there exists a
Hilbert space $\hil_0$ with dimension $m-d$ and an isometric
isomorphism $U:\FF^m\rightarrow \hil\oplus \hil_0$ such that
\begin{equation}\label{el U}
UGU^*=\begin{pmatrix} S & 0\\ 0& 0
\end{pmatrix}\hspace{-0.2 cm}\barr{c} \hil\\\hil_0 \earr
\end{equation}
 Therefore, $(\|\phi_i\|^2)_{i=1}^m
 \prec (\sigma(S^\cF),0_\sim)$
 where $0_\sim \in \RR^{m-d}$. 
\end{pro}
As a consequence of Proposition \ref{el s y el g} we see that, if $\sigma(G)\in \FF^m$ (resp $\sigma(S)\in \CC^d$) denote the eigenvalues of $G$ counted with multiplicity
then $\sigma(G)=(\sigma(S),0_\sim)$ where $0_\sim\in \FF^{m-d}$.

\begin{teo}[\cite{JDM,MR}]\label{hay frame} Let $S\in \mathcal M_d(\FF)^+$ and let
$\a=(a_i)_{i=1}^m$ be a sequence of positive numbers.
Then, there exists a  sequence $\{\phi_i\}_{i=1}^m\subset \hil$ with frame operator $S$ and such that $\|\phi_i\|=a_i$
for every $1\leq i\leq m$ if and only if $$\sum_{i=1}^k a_i^2\leq
\sum_{i=1}^k \lambda(S)_i,\ \text{ for } \ 1\leq i\leq d-1, \
\text{ and }\ \sum_{i=1}^m a_i^2=\tr(S).$$
\end{teo}

\section{Convex functions defined on frame operators.}\label{Convex funct}
In this section we define a family  functions $P_f$ on the set of frame operators of sequences $\cF$ in $\CC^d$, starting from a convex function $f:\RR_{\geqslant 0}\rightarrow \RR_{\geqslant 0}$. As a particular case, we recover the frame potential, introduced by Benedetto and Fickus, in \cite{BF} with a specific convex function $f$.

When we restrict our attention to special sets of  sequences, namely, those sequences with a prescribed set of norms, we are able to compute the minimum value taken by $P_f$ on the  corresponding set of frame operators and to characterize the spectrum of minimizers of  $P_f$, for every $f$ non decreasing and convex function which  satisfies $f(0)=0$.

\begin{fed}\label{defi de potencial}
Let $f:\RR_{\geqslant 0}\rightarrow \RR_{\geqslant 0}$ be a non decreasing convex function.
Then, the frame potential associated to $f$,
denoted $P_f$, is the functional defined on the set of frame operators of 
sequences  in $ \CC^d$ given by
\begin{equation}\label{potencial}
P_f(S^\cF)=\tr(f(S^\cF))
\end{equation}for every $\cF=\{\phi_i\}_{i=1}^m\subset \CC^d$. In detail, if we denote by $\lambda =(\lambda_i)_{i=1}^d$ the eigenvalues of $S^\cF$ counted with multiplicity, then
$P_f(S^\cF)=\sum_{i=1}^d f(\lambda_i)$.
\end{fed}

\begin{rem}\label{la gram pot}
Using the relation between $G^\cF$ and $S^\cF$ shown in Proposition \ref{el s y el g}, we have
$$  \tr(f(G^{\cF}))=P_f(S^\cF)+(m -d)f(0)$$
In particular, if $f(0)=0$, $P_f(S^\cF)$ can be computed using the Grammian matrix.
\end{rem}

\begin{exa}[Benedetto-Fickus's potential]\label{potencial usual}
Let $f:\RR_{\geqslant 0}\rightarrow \RR_{\geqslant 0}$ be the strictly convex function
$f(x)=x^2$. Then, the frame potential associated to $f$ is
\begin{equation}\label{benedetto-fickus}P_f(S^\cF)=\tr((S^\cF)^2)=\tr((G^\cF)^2)=\sum_{i,\,j=1}^ m|\langle
\phi_i,\phi_j\rangle|^2\end{equation} that is, the frame potential
as defined by Benedetto and Fickus in \cite{BF}
\end{exa}

In what follows, given $\cF=\{\phi_i\}_{i=1}^m\subset \CC^d$ and $\alpha\in \CC$ we denote by $\alpha \cF=\{\alpha\,\phi_i\}_{i=1}^m$. On the other hand, given $\cF_1=\{\phi_i\}_{i=1}^{M_1}$, $\cF_2=\{\psi_i\}_{i=1}^{M_2}\subset \CC^d$ then $\cF_1\sqcup\cF_2$ denotes the list of $M_1+M_2$ vectors obtained by juxtaposition of $\cF_1$ and $\cF_2$. Note that, if $\cG=\alpha \cF_1\sqcup \cF_2$ then $$S^\cG=|\alpha|^2\,S^{\cF_1}+ S^{\cF_2}.$$
\begin{teo}
\label{algunas implicaciones}
Let $f:\RR_{\geqslant 0}\rightarrow\RR$ be a non decreasing convex function  and $\cF_1=\{\phi_i\}_{i=1}^{M_1}$, $\cF_2=\{\psi_i \}_{i=1}^{M_2} \subset \CC^d$. 
\begin{enumerate}
\item  If $S^{\cF_1} \prec S^{\cF_2}$ then $$P_f(S^{\cF_1})\leq
P_f(S^{\cF_2}).$$
\item Assume further that $f$ is a strictly convex function, $S^{\cF_1} \prec S^{\cF_2}$
and $P_f(S^{\cF_1})$ = $ P_f(S^{\cF_2})$.
Then, there exists a unitary operator $U\in \cM_d(\CC)$ such that
$$US^{\cF_1}U^*=S^{\cF_2}.$$
\item If $t\in [0,1]$ and $\cG=t^{1/2} \cF_1\sqcup (1-t)^{1/2} \cF_2$ then  
$$P_f(S^\cG)\leq tP_f(S^{\cF_1})+ (1-t)P_f(S^{\cF_2}).$$
\item If $\cG= \cF_1\sqcup  \cF_2$ then  $$P_f(S^\cG)\geq P_f(S^{\cF_1})+P_f(S^{\cF_2}).$$
\end{enumerate}
\end{teo}
\begin{proof}
The first two items are well known (see \cite{Ba,{horn2}}). The last two inequalities above are also well known (see \cite[Theorem 1-24]{AU}) for these functionals. 
\end{proof}

\begin{rem}
For $g=-f$,  $P_g(S)=\tr (g(S))$ for $S\in \cM_d(\CC)^+$ are  called ``entropy-like"  functionals in \cite{AU}. Notice that the minimization of the functions $P_f$ corresponds to the maximization of the entropy-like functional $P_g$.
\end{rem}

Let $c>0$ and $\a=\{a_i\}_{i=1}^m$ be a sequence of positive elements arranged in decreasing order.  In what follows we shall consider the following sets:
$$\cA(c)=\{\cF=\{\phi_i\}_{i=1}^m\subset \CC^d,\,\sum_{i=1}^m\|\phi_i\|^2=\tr(S^\cF)=c\},$$
$$\cB(\a)=\{\cF=\{\phi_i\}_{i=1}^m\subset \CC^d,\, \|\phi_i\|^2=a_i \text{ for every } i\}.$$
Observe that, by Theorem \ref{hay frame}, the sets of frame operators for  sequences in $\cA(c)$ and $\cB(\a)$ can be well characterized:
\begin{equation}\label{el te}
\cT(c)=\{ S^\cF, \, \cF\in \cA(c)\}=\{ S\in \cM_d(\CC)^+, \, \la(S)\in \K(c)\}.
\end{equation}
\begin{equation}\label{el ere}
\cR(\a)=\{ S^\cF, \, \cF\in \cB(\a)\}=\{ S\in \cM_d(\CC)^+, \, \la(S)\in \cP(\a)\}.
\end{equation}

\begin{teo}\label{sin normas}
Let  $f:\RR_{\geqslant 0}\rightarrow \RR_{\geqslant 0}$ be a non decreasing convex function and $P_f$ the functional associated to $f$ and let $c>0$.
Then, if $\cF \in \cA(c)$ is a tight frame, then 
\[ P_f(S^\cF)\leq P_f(S^\cG) \quad \forall S^\cG \in \cT(c). \]
Moreover, if in addition $f$ is strictly convex and  $S^\cF$ is a local minimum of $P_f$ considering the operator norm in $\cT(c)$, then $S^\cF=\frac{c}{d} I$  so $\cF$ is a tight frame.
\end{teo}

\begin{proof}

The proof follows immediately from  Proposition \ref{maximos y minimos}. Indeed, $S^\cF\in \cT(c)$ is a global minimum for $P_f$ if and only if $\la(S^\cF)=\v$, i.e. $S^\cF=\frac{c}{d}I$, which means that $\cF$ is a tight frame in $\cA(c)$. On the other side, if $\la=\la(S^\cF)\neq \v$, then by Prop. \ref{maximos y minimos} for every $\eps>0$ sufficiently small, there exist $\la_\eps \in \K(c)$ such that $\la_\eps\prec \la$, $\lambda^\downarrow\neq \lambda_{\eps}^\downarrow$ and $\|\la -\la_\eps\|<\eps$. Thus, if  $S^\cF=U^*\diag (\la )U$ with $U$ unitary,  it is clear that $S_\eps=U^*\diag(\la_\eps )U\in \cT(c)$ satisfies $\|S^\cF- S_\eps\|<\eps$ and $P_f(S_\eps)<P_f(S^\cF)$, by Thm. \ref{algunas implicaciones}.

\end{proof}

\begin{teo}\label{T:minimos globales}
Let  $f:\RR_{\geqslant 0}\rightarrow \RR_{\geqslant 0}$ be a non decreasing convex function and $P_f$ the functional associated to $f$. Let $\a =(a_i)_{i=1}^m$ be a non increasing sequence of strictly positive numbers with $d\leq m$.  
Suppose that $\cF \in \cB(\a)$ is of the form
\begin{equation}\label{minimos globales}
 \{ \sqrt{a_i} \,
e_i\}_{i=1}^r\cup \{\phi_i\}_{i=r+1}^m
\end{equation}
where $\{e_i\}_{i=1}^d$
is an o.n.  basis for $\CC^d$, $r$ is the
\emph{d-irregularity} of $\a$ and
  $\{\phi_i\}_{i=r+1}^m$ is a
tight frame for $\text{span}\{e_i\}_{i=r+1}^d$ with frame constant $\displaystyle{c=(d-r)^{-1}\sum_{i\geq r+1}a_i}$.

\smallskip

Then, $S^\cF$ is a global minimum for $P_f$ in $\cR(\a)$. Moreover, if $f$is strictly convex and  $S^\cF$ is a local minimum for $P_f$ in $\cR(\a)$ (considering the operator norm), then $\cF$ is as in \eqref{minimos globales}.

\end{teo}

\begin{proof}
Let $\cF \in \cB(\a)$ be of the form given in \eqref{minimos globales}. Therefore, the (ordered)  spectrum of the frame operator $S^{\cF}$ is $\v=(a_1,\ldots,a_r, c, \ldots, c)$ where $c=(d-r)^{-1}\sum_{i\geq r+1}a_i$ is an eigenvalue with multiplicity $d-r$.
Then,  by the  Proposition \ref{hay min para la mayo} and Theorem \ref{algunas implicaciones}, we can conclude that $\cF$ is a global minimum in $\cB(\a)$.

Now, let $\cG=\{\psi_i\}_{i=1}^m\in \cB(\a)$ be such that $\la (S^\cG)=\v$. Then the (optimal) upper frame bound of $\cG$ is $a_1$ and we have  
\[
\norm{\psi_1}^4+\sum_{j>1 }\left|\pint{ \psi_j \, , \,  \psi_1}\right|^2 \leq a_1 \norm{\psi_1}^2=\norm{\psi_1}^4
\]
Therefore, $\psi_1$ is orthogonal to $\psi_j$ for $j\neq 1$. By restriction to $\sspan \{\psi_i\}_{i=2}^m$, we deduce that $\pint{\psi_2\ , \ \psi_i}=0$ for $i\neq 2$ in the same way. Therefore  we can conclude that $\pint{\psi_i\ , \psi_j}=0$ for every $1\leq i\leq r$, $j\neq i$, in particular we  define the orthonormal set $e_i=a_i^{-1/2}\psi_i$ for $1\leq i\leq r$. We then complete it to an o.n.b. $\{e_i\}_{i=1}^d$.

Finally, since the frame operator restricted in the orthogonal complement of the space spanned by  $\{e_i\}_{i=1}^r$ is a multiple of the identity, the rest of the frame is a tight frame in its span. Then, $\cG$ can be described as in \eqref{minimos globales}. 

Let $S^\cF\in \cR(\a)$ be such that $\la(S^\cF)$ is not $\v\in \cP(\a)$. Therefore, by the last statement of  Prop. \ref{hay min para la mayo} and arguing as in Thm. \ref{sin normas}, given $\eps>0$, we can find a positive definite operator $S_\eps\in \cR(\a)$ such that $\|S_\eps-S^\cF\|<\eps$, $\lambda(S_\eps)\neq \lambda(S^\cF)$ and $S_\eps \prec S^\cF$. Then $P_f(S_\eps)<P_f(\cF)$ for every strictly convex function $f$, by Theorem \ref{algunas implicaciones}. In particular, by the previous paragraph, every local minimum for $P_f$ in $\cR(\a)$ is a global minimum, so it is a frame operator of a frame given by \eqref{minimos globales}.  

\end{proof}

Theorem (A) in the Introduction is now a consequence of the identities \eqref{el te}, \eqref{el ere} and Theorems \ref{sin normas}, \ref{T:minimos globales}.

\begin{cor}\label{minimizadores de potenciales}
Let $\cF=\{\phi_i\}_{i=1}^m\in \cA(c)$ and let
$f:\RR_{\geqslant 0}\rightarrow \RR_{\geqslant 0}$ be a non decreasing convex function.
We have the following inequalities:
\begin{equation}\label{primera cota} (d-1)\cdot
f(0)+f(c)\geq P_f(S^\cF)\geq d\cdot f(\frac{c}{d}),
\end{equation}

And, for $\cF \in \cB(\a)$ we have
\begin{equation}\label{tercera cota}
(d-1)\cdot
f(0)+f(\sum_{i=1}^m a_i)\geq P_f(S^\cF)\geq \sum_{i=1}^r f(a_i)+ (d-r)\cdot f(h)
\end{equation}
with $h=(d-r)^{-1}\sum_{i\geq r+1} a_i$. Moreover, if in addition $f$ is strictly convex and the lower bound is attained in \eqref{primera cota} (respectively in \eqref{tercera cota}) then $\cF$ is a tight frame (respectively is as in \eqref{minimos globales} for some o.n.b. $\{ e_i\}_{i=1}^d$ for $\CC^d$).
\end{cor}

\subsection{Some applications of the previous results}

Let us begin with the following example in order to illustrate the
content of our previous results.

\begin{exa}[continuation of example \ref{potencial usual}] Let
$f(x)=x^2$ and note that, by equation (\ref{benedetto-fickus}), if
$\cF=\{\phi_i\}_{i=1}^m\in \cA(c)$ then
$$P_f(S^\cF)=\sum_{i,j=1}^m |\langle
\phi_i,\phi_j\rangle|^2.$$  
Note that $f$ is a
strictly convex function and that $f(\lambda\cdot
x)=\lambda^2\cdot f(x)$ for every $\lambda\geq 0$, so we can take
$g(\lambda)=\lambda^2$. Then, equation (\ref{tercera cota})
becomes
\begin{equation}\label{la ecua de Wal} 1\geq
\frac{\sum_{i=1}^m|\langle
\phi_i,\phi_j\rangle|^2}{(\sum_{i=1}^m\|\phi_i\|^2)^2}\geq d\cdot
\frac{1}{d^2}=\frac{1}{d} \end{equation} which is the
\emph{generalized Welch inequality} of \cite{Wal}. Moreover, by
Theorem \ref{minimizadores de potenciales} we deduce that the
lower bound (resp the upper bound) in equation (\ref{la ecua de
Wal}) is attained if and only if $\cF$ is a tight frame with frame
bound $\frac{c}{d}$ (resp if and only if $\lin(\cF)$ has dimension
1).

\end{exa}

Of course, the function $f(x)=x^2$ is probably the most simple
function that can be used to produce a reasonable frame potential.
 In the following examples we shall
investigate other choices of convex functions.

\begin{exa}[$n$-th frame potential] Let $n\geq 2$ and consider
$f_n(x)=x^n$ for $x\geq 0$. Then, $f$ is an increasing strictly
convex function and produce the \emph{$n$-th frame potential}
given by $$P_n(S^\cF)=\tr((S^\cF)^n)$$ where $S^\cF$ is
the frame operator of the sequence $\cF= \{\phi_i\}_{i\in m}\subset
\CC^d$. Since $f(0)=0$ then
 we
have
\begin{equation}\label{el potencial
n}P_n(\{\phi_i\}_{i=1}^m)=\tr((G^\cF)^n)=\sum_{i_1,\ldots,\,i_n=1}^m
\ \prod_{j=1}^{n}\langle \phi_{i_j},\phi_{i_{j+1}}\rangle
\end{equation}where we follow the convention $i_{n+1}=i_1$. Note
that $P_2$ is the usual frame potential. Indeed, formula (\ref{el
potencial n}) is a consequence of the identity
\begin{equation}\label{la identidad}\langle (G^\cF)^ne_k,e_k\rangle=
\sum_{\begin{subarray}{c} i_2,\ldots,\,i_n=1\\ i_1=k \end{subarray}}^m \ \prod_{j=1}^{n}\langle
\phi_{i_j},\phi_{i_{j+1}}\rangle\geq 0
\end{equation}
In this case, using equation (\ref{el potencial n}), equation
(\ref{tercera cota}) becomes
\begin{equation}\label{otra cota}
1\geq \frac{\sum_{i_1,\ldots,\,i_n=1}^m \ \prod_{j=1}^{n}\langle
\phi_{i_j},\phi_{i_{j+1}}\rangle}{(\sum_{i=1}^m\|\phi_i\|^2)^n}\geq
\frac{1}{d^{n-1}}
\end{equation}while equation (\ref{la identidad}) implies
\begin{equation}\label{mas desigualdades}
\max_{1\leq k\leq m} \sum_{\begin{subarray}{c} i_2,\ldots,\,i_n=1\\ i_1=k \end{subarray}}^m \
\prod_{j=1}^{n}\langle \phi_{i_j},\phi_{i_{j+1}}\rangle\geq
\frac{(\sum_{i=1}^m\|\phi_i\|^2)^n}{m\cdot d}
\end{equation}
As before, the lower bound in formula (\ref{otra cota}) is
attained if and only if $\cF$ is a tight frame with frame bound
$\frac{c}{d}$. Analogously, the bound in equation (\ref{mas
desigualdades}) is attained if and only if $\cF$ is a tight frame.

\end{exa}

\begin{exa}[von Neumann  Entropy]
If we consider the concave function $f(x)=-x\ln(x)$, then $P_f$ restricted to density matrices is the well known von Neumann entropy in quantum information theory. Roughly speaking, it measures the lack of information about the state of a system. Theorems \ref{sin normas} and \ref{T:minimos globales} show, as a particular case, the structure of maximizers of the entropy without restrictions in the first case and with the restriction: $\{ S\mbox{ a density matrix with  } (\la(S),0_{m-d}) \succ \a \}$ for a fixed positive sequence $\a$ with $\sum_{i=1}^m a_i=1$.  
\end{exa}

\subsection{Convex functions over CGU frames.}
In this section we use the previous techniques to characterize the global minimizers of $P_f$ when restricted to the compound geometrically uniform frames, with a prescribed list of norms.

 \begin{fed} 
Let $\mathsf{G}$ be a finite abelian group of unitaries in $\cM_d(\CC)$, and $\varphi\in \CC^d$.  If the set 
$ \mathsf{G}\cdot \varphi=\{ U\varphi : U\in  \mathsf{G}\}$ is a frame the we say that $ \mathsf{G}\cdot \varphi$ is a {\bf geometrically uniform frame} (GU). When $\mathsf{G}$ acts on a larger set of functions, $\Phi=\{\varphi_i \in \CC^d : 1\leq i\leq m\}$ and $\mathsf{G}\cdot \Phi$ is a frame, we say that it is a {\bf compound geometrically uniform frame} (CGU).
\end{fed}

From now on, in order to simplify the computations, we assume also that $ \mathsf{G}$ is cyclic. Let suppose then that we have $ \mathsf{G}=\{U^i : 0\leq i\leq n-1 \}$, where $U$ is a unitary such that $U^n=I$. Thus, we shall consider frame sequences of the form $\cF= \mathsf{G}\cdot\Phi=\{U^i\varphi_j: 0\leq i\leq n-1 \, , \quad 1\leq j\leq m\}$.

We are interested in minimizing $P_f$ when we restrict $P_f$ to the set of frame operators of CGU frames:

\[  \mathsf{G} \cdot \cB(\a)=\{ \mathsf{G} \cdot \cF \, : \, \cF\in \cB(\a)\},\] 
where $ \mathsf{G}$ is a fixed cyclic group of unitaries, $\a$ is fixed.
Clearly $ \mathsf{G} \cdot \cB(\a)\subset \cB(\b)$, where $\b =\{b_i\}_{i=1}^{nm}$ is the sequence $\a$ repeated $n$ times. Then, by Corollary \ref{minimizadores de potenciales}, if $\cF\in  \mathsf{G}\cdot \cB(\a)$,
 \begin{equation}\label{repetidos} 
P_f(S^\cF)\geq \sum_{i=1}^r f(b_i)+ (d-r)\cdot f(h),
\end{equation} 
where $h=(d-r)^{-1}\sum_{i=r+1}^{nm}b_i$ and $r$ is the $d$-irregularity of $\b$.
The previous inequality can be stated in terms of $\a$ if we characterize the $d$-irregularity of $\b$.
\begin{pro}\label{irregularidad del b}
Let $\a=(a_i)_{i=1}^m$ be a non increasing sequence of positive numbers and let $\b=(b_i)_{i=1}^{nm}$ be a sequence given by:
\[b_j=a_i \quad \text{ for } \quad j=(i-1)n+s, \quad 1\leq s\leq n, \; 1\leq i\leq m.\]
Then, if $r_0$ is the $d$-irregularity of $\b$, $r_0=nr$, where 
\[r=\max \{j: (\frac{d}{n}-j)\,a_j>\sum_{k=j+1}^m a_k \}. \] 
\end{pro}

\begin{proof}
The result is clear if $r_0=0$. If $r_0\neq 0$, then it holds that $n$ divides $r_0$. Indeed, by definition of $r_0$, $b_{r_0}\neq b_{r_0+1}$ which can only occur if $r_0=nr$, $r\in \{1,\ldots m\}$.
Finally, 
$$ r_0=\max \{nj : \, (d-nj)\,b_{nj}> \sum_{k=nj+1}^{nm}b_k\} = n\,\max\{ j : \, (\frac{d}{n} -j)\,a_j >\sum_{k=j+1}^m a_k\}.$$
\end{proof}

\begin{teo}\label{minimos en CGU}
Let $ \mathsf{G}$, $\a$ and $\cB(\a)$ as before. Suppose that $n|d$   and that there exists an orthonormal family $\{e_i\}_{i=1}^{N}$, with $N=\frac{d}{n}$ such that the set $\{U^ke_j \; 1\leq k\leq n, \, 1\leq j\leq N \}$ is an orthonormal basis of $\C^d$. 
Let $\cF$ in $\cB(\a)$ be of the form
\begin{equation}\label{globales en cgu}
\cF'= \{ \sqrt{a_i} \,
b_i\}_{i=1}^r\cup \cD
\end{equation}
where $\cE=\{b_i\}_{i=1}^r$
is an orthonormal set such that $ \mathsf{G}\cdot \cE$ is orthonormal, $r$ is the
\emph{$N$-irregularity} of $\a$ and
  $ \mathsf{G} \cdot\cD$ is a
tight frame for $\ds{\sspan\left(  \mathsf{G}\cdot \cE \right)^\perp}$  with frame constant $h=(N-r)^{-1}\sum_{k=r+1}^m a_k$. Denote $\cF= \mathsf{G}\cdot \cF' \in  \mathsf{G}\cdot\cB(\a)$.
\medskip

Then $S^\cF$ is a global minimum for $P_f$ in the set of frame operators of $ \mathsf{G} \cdot \cB(\a)$. Conversely, if in addition $f$ is strictly convex, and $S^\cF$ is a global minimum for $P_f$, then $\cF$ is of the form $ \mathsf{G} \cdot \cV$, with $\cV$ as in \eqref{globales en cgu}. 

\end{teo}

\begin{proof}
By Thm. \ref{T:minimos globales} and Prop. \ref{irregularidad del b} it is clear that if such  sequence exists, then $S^\cF$ is a global minimum in $\cR(\b)$ (using the previous notation), so it is a global minimum when we restrict $P_f$ to the frame operators of  $ \mathsf{G} \cdot \cB(\a)$. Moreover, if $f$ is strictly convex, every global minimum must be of this form, by Thm. \ref{T:minimos globales}.

Then, in order to prove the statement we need to show that such  sequence exists.  Indeed let $\cF'$ be the sequence given by
\[\{ \sqrt{a_i} \,
e_i\}_{i=1}^r\cup \{\phi_i\}_{i=r+1}^m
\]
where $\{e_i\}_{i=1}^N$
is the orthonormal set existing by the hypotheses, $r$ is the
\emph{$N$-irregularity} of $\a$ and
  $\{\phi_i\}_{i=r+1}^m$ is a
tight frame for $\sspan \{ e_k\}_{k=r+1}^N$, with frame constant $h=(N-r)^{-1}\sum_{k=r+1}^m a_k$. Such frame exists  by Theorem \ref{T:minimos globales}.

Clearly, for every $1\leq k\leq n$, the set $\{U^k\phi_i\}_{i=r+1}^m$ is a tight frame (with the same constant $h=(N-r)^{-1}\sum_{k=r+1}^m a_k$) for $\sspan \{U^ke_i\}_{i=r+1}^N$, therefore, $\cD= \mathsf{G}\cdot \{\phi\}_{i=r+1}^m$ is a tight frame of $\ds{\sspan\left(  \mathsf{G}\cdot \cE \right)^\perp}$  with frame constant $c=(N-r)^{-1}\sum_{k=r+1}^m a_k$, where $\cE=\{e_i\}_{i=1}^r$.

\end{proof}
\begin{rem}
If in addition we assume that the initial vectors $\cF$ lie on the $\frac{d}{n}$- dimensional subspace $\cK$ generated by $\{e_i\}_{i=1}^N$ ($\frac{d}{n}=N$) of $\CC^d$ we can conclude that the global minimizers are of the form given in \eqref{globales en cgu}, where $r$ is the $\frac{d}{n}$-irregularity of $\a$ and $\cD$ forms a tight frame on $\cK\cap (\sspan \{b_i\}_{i=1}^r)^\perp$. Indeed, in this case the Grammian matrix of $\mathsf G\cdot \cF$ is block-diagonal.

A special case of this situation is given on convolutional frames studied in \cite{FJKO}. In particular, previous Theorem can be seen as a partial generalization to \cite[Thm. 6]{FJKO}.
 
\end{rem}

\begin{cor}\label{minimizadores de potenciales en CGU}
Under the hypotheses of Theorem \ref{minimos en CGU}, for  $\cF \in \cB(\a)$ we have

\begin{equation}\label{tercera cota CGU}
(d-1)\cdot
f(0)+f(n \cdot \sum_{i=1}^m a_i)\geq P_f(S^{\mathsf{G}\cdot \cF})\geq n \left\{ \sum_{i=1}^r f(a_i)+ (d-r)\cdot f(h)\right\}
\end{equation}
with $h=(d-r)^{-1}\sum_{i\geq r+1} a_i$. Moreover, if in addition $f$ is strictly convex and the lower bound is attained in  \eqref{tercera cota CGU} then $\cF$ is as in \eqref{globales en cgu}.
\end{cor}

\section{ From frame operators to frames.}\label{frame op to frames}

In the previous section we have considered the function $P_f$ associated to a convex function $f$ as a function of the frame operators; we have described the structure of local minimizers of $P_f$ when restricted to the sets $\cT(c)$ and $\cR(\a)$ with respect to the norm topology. 

We are now interested in considering $P_f$ defined on frames  
\[ P_f(\cF):=P_f(S^\cF)=\tr(f(S^\cF))\]
for $\cF=\{\phi_i\}_{i=1}^m\subset \CC^d$, and studying the structure of global and local minimizers of these functions when restricted to the sets $\cA(c)$ and $\cB(\a)$, with respect to the vector-vector distance \begin{equation}d(\cF,\cG)=\max_{1\leq i\leq m}\|\phi_i-\psi_i\|\end{equation} for  sequences $\cF=\{ \phi_i \}_{i=1}^m$, $\cG=\{ \psi_i \}_{i=1}^m$. It is worth noting that 
the norm distance between frame operators can not bound the vector-vector distance; indeed if $\sigma$ is a permutation of order $m$ and $\cG=\{f_{\sigma(i)}\}_{i=1}^m$ then $S^\cF=S^\cG$ while $d(\cF,\cG)\neq 0$ possibly. This implies that the results in the previous section can not be used to obtain a complete characterization of the local minimizers in this new setting.

Our approach to this new point of view involves the study of the existence local cross sections of the map $\cF\mapsto S^\cF$
when it is restricted to  $\cA(c)$ and $\cB(\a)$ respectively (note that the restriction on the norms which defines $\cB(\a)$ is a condition on the main diagonal of $G^\cF$).

To begin with, Theorem \ref{sin normas} implies that if a  sequence $\cF=\{\phi_i\}_{i=1}^m\in \cA(c)$ does not have the structure of a local (global) minimizer of $P_f$ on $\cT(c)$, for a strictly convex function $f$, then for every $\eps>0$ there exists a $S\in \cT(c)$ such that $\|S-S^\cF\|\leq \eps$ and $P_f(S)<P_f(S^\cF)$. In order to show that $\cF$ is not a local minimum of $P_f$ on $\cA(c)$ with respect to the vector-vector distance the following problem arises: given such $S$, is there any  sequence $\cG=\{\psi_i\}_{i=1}^m\in\cA(c)$ such that $S^G=S$ and $d(\cF,\cG)\leq \delta(\eps)$, with $\lim_{\eps\rightarrow 0}\delta(\eps)=0$?
A positive answer to this question is given in the following

\begin{pro}\label{pert f op sin normas}
Let $\cF=\{\phi_i\}_{i=1}^m \in \cA(c)$ and let $S\in \cT(c)$ be such that $\|S-S^\cF\|<\eps$. Then there exist $\cG=\{\psi_i\}_{i=1}^m\in \cA(c)$ such that $d(\cF,\cG)<\eps \rai$ and $S^\cG=S$.
\end{pro}
\begin{proof}
Consider $T^\cF=(S^\cF)\rai W$ the polar decomposition of $T^\cF$. Then, since $\|S-S^\cF\|<\eps \, $,  $\, \|S\rai -(S^\cF)\rai\|<\eps \rai$ by \cite[Thm. X.1.1]{Ba}).

Now let $\cG=\{\psi_i\}_{i=1}^m$, where $\psi_i=S\rai W e_i$ for $1\leq i\leq m$. Then $T^\cG=S\rai W$, $S^\cG=T^\cG (T^\cG) ^*=S$ and for $1\leq i\leq m$
$$\|\psi_i-\phi_i\|\leq \|T^\cG-T^\cF\|\leq \|S\rai -(S^\cF)\rai\|<\eps \rai.$$
\end{proof}
The previous result combined with Theorem \ref{sin normas} provide a complete characterization of the local (global) minimizers of $P_f$ on $\cA(c)$ with respect to the vector-vector distance, for a strictly convex $f$.

\begin{teo}\label{el sin label}
Let $f: \RR_{\geqslant 0}\rightarrow \RR_{\geqslant 0}$ be a non decreasing convex function. If $\cF\in \cA(c)$ is a tight frame then it is a global minimizer of $P_f$ on $\cA(c)$. Moreover, if $f$ is a strictly convex function then every local minimum of $P_f$ on $\cA(c)$ with respect to the vector-vector distance is a tight frame. 
\end{teo}
\begin{proof}
The first part of the statement follows from Theorem \ref{sin normas} and \eqref{el te}.
By the proof of \ref{sin normas}, if $\cF\in \cA(c)$ is not tight, then  for every $\eps>0$, there exists $S_\eps\in \cT(c)$  such that $\|S^\cF-S_\eps\|<\eps^2$ and $P_f(S_\eps)<P_f(S^\cF)$. 
Finally, by Proposition \ref{pert f op sin normas}, there exist $\cG=\{\psi_i\}_{i=1}^m\in \cA(c)$ such that $S^\cG=S_\eps$ and $\|\phi_i-\psi_i\|<\eps$. 
\end{proof}

As before, in order to obtain a characterization of local minimizers of $P_f$ on $\cB(\a)$ with respect to the vector-vector distance using Theorem \ref{T:minimos globales} we are led to consider the following perturbation problem: given a  sequence $\cF=\{\phi_i\}_{i=1}^m\in\cB(\a)$ and $S\in \cR(\a)$ with $\|S^\cF-S\|\leq \eps$, is there a  sequence $\cG=\{\psi_i\}_{i=1}^m\in\cB(\a)$ with $S^\cG=S$ and $d(\cF,\cG)\leq \delta(\eps)$ with $\lim_{\eps\rightarrow 0}\delta(\eps)=0$. The constrain $\cG\in\cB(\a)$ seems to be hard to deal with. For example, notice that we have no control on the norms of the vectors in $\cG$ constructed in Proposition \ref{pert f op sin normas}. 
On the other hand, it is convenient to work with the Grammian since the restriction $\cG\in \cB(\a)$ is equivalent to $\mathrm d(G^\cG)=\a$, where $\mathrm d(X)\in \CC^m$ 
denotes the main diagonal of the $m\times m$ complex matrix $X$. 

We have only obtained  partial results which are presented in the following Proposition. The proof  depends strongly on geometrical aspects and it is developed in the appendix. 

\begin{pro}\label{prop: anda con sec loc}
Let $\cF=\{\phi_j\}_{j=1}^m\subseteq \CC^d$ be a frame, let $S=S^\cF$ be its frame operator and assume that $\cF$ can not be partitioned in two sets of mutually orthogonal vectors. Let $\{S_i\}_{i}\subseteq \mathcal M_d(\CC)^+$ be a sequence converging to $S$. Then, for every $\eta>0$ there exists $i_1\in \NN$ such that for each $i\geq i_1$ there exists a frame $\cG(i)=\cG=\{ \psi_j\}_{j=1}^m$ such that: 

\begin{enumerate}
\item
 $\|\psi_j\|=\|\phi_j\|$ for $1\leq j\leq m$.
\item
$\|\psi_j-\phi_j\| \leq \eta$ for $1\leq j\leq m$.
\item
$S^\cG=S_i$.
\end{enumerate}
\end{pro}

\begin{teo}\label{estruc min loc caso gen}
Let $f: \RR_{\geqslant 0}\rightarrow \RR_{\geqslant 0}$ be a non decreasing convex function. If $\cF\in \cB(\a)$ has the structure as in \eqref{minimos globales} then it is a global minimizer of $P_f$ on $\cB(\a)$. 

 If in addition $f$ is a strictly convex function, then every global minimum of $P_f$ on $\cB(\a)$ is as in \eqref{minimos globales}. 
Moreover, for such $f$ then every $\cF=\{\phi_i\}_{i=1}^m\in \cB(a)$ such that it can not be partitioned in two mutually orthogonal sets of vectors is a local minimum if and only if is a global minimum. 
\end{teo}  

\begin{proof}
The first part of the statement follows from Theorem \ref{T:minimos globales} and \eqref{el ere}. 

Assume now that $\cF$ is not a global minimum; by the proof of Thm. \ref{T:minimos globales}, there is a sequence of operators $\{S_n$\} such that $S_n$ converges to $S^\cF$ and such that $P_f(S_n)<P_f(S^\cF), \forall n$.

Let $\eps>0$, then, by Thm. \ref{anda con sec loc}, for a sufficient large $n_0\in \NN$, there exist a frame $\cG=\{\psi_i\}_{i=1}^m\in \cB(\a)$ such that $\|\phi_i-\psi_i\|<\eps$ and $S^\cG=S_n$. In particular, $P_f(\cG)<P_f(\cF)$.  
\end{proof}
Theorem (B) in the Introduction follows immediately from Theorems \ref{el sin label} and \ref{estruc min loc caso gen}. 

It is clear that $P_f(\cF)=P_f(\cF_1)+P_f(\cF_2)$ if $\cF=\cF_1\cup \cF_2$ with the vectors in $\cF_2$ being orthogonal to those in $\cF_1$ (we shall denote this by $\cF_1\perp \cF_2$).  This simple observation and the  previous result, allows a  reduction of the set of possible local minimizers for $P_f$:

\begin{cor}\label{cor estr min loc caso gen}
Let $\cF\in \cB(\a)$ such that $\cF=\cF_1\cup \cF_2$ with $\cF_1=\{\phi_i\}_{i=1}^{M_1}\perp \cF_2=\{\phi'_i\}_{i=1}^{M_2}$, and suppose that $\cF_1$ can not be partitioned into two mutually orthogonal sequences and it is not a global minimizer for $P_f$ restricted to the set
\[ \cB(\a_1)=\{\{\psi_j\}_{j=1}^{M_1} \, : \, \psi_j\in \sspan \cF_1, \, \|\psi_j\|=\|\phi_{j}\| \quad  1\leq j\leq M_1\}.\] 
Then, $\cF$ is not a local minimizer for $P_f$.
\end{cor}

Note that the general structure of local minimums of arbitrary function $P_f$ can not be inferred from Theorem \ref{estruc min loc caso gen} and Corollary \ref{cor estr min loc caso gen}. Still, these results allow to a reduction of the general situation to a particular case (see Problem ($\star$) below). In order to exemplify the ideas involved, we recover \cite[Theorem 10]{casazza2} about the structure of general minimizers in the particular caso of the Benedetto-Fickus potential.

\begin{teo}\label{teo:recovering} 
Any local minimizer $\cF=\{\phi_i\}_{i=1}^m$ of the Benedetto-Fickus potential in $\cB(\a)$ with respect to the distance $d(\cdot , \cdot)$ is a global minimizer of this potential and hence has the structure given in \eqref{minimos globales}
\end{teo}

\begin{proof}
Suppose that we have a frame $\cF\in \cB(\a)$ which is not a global minimum for the Benedetto - Fickus potential  FP. We must show that then it is not a local minimum.

Let $\cF=\cF_1\cup \ldots \cup \cF_k$ its minimal  decomposition in pairwise orthogonal subsets (minimal in the sense that no $\cF_j$ can be partitioned in two mutually orthogonal subsets). By Corollary \ref{cor estr min loc caso gen}, if there exist $1\leq i\leq k$ such that $\cF_i$ is not a global minimum for FP (restricted to $\cB(\a_i)$), then $\cF$ is not a local minimizer. 

So  we can suppose that every  $\cF_i$ is a global minimum on $\cB(\a_i)$. Then by Theorem \ref{T:minimos globales},  $\cF_i$ is   tight on its span (possibly with a single vector), with frame constant $c_i$, for $1\leq i\leq k$. We claim that in this case, there is a pair $i,j$ such that the frame $\cF_i\cup \cF_j$ is not a global minimum for FP on span $\cF_1\cup \cF_2$ with the restriction given by the vector norms in $\cF_i$, $i=1,2$.

Indeed, if there exists a pair $\cF_i$ and $\cF_j$, each with  two or more vectors, and with constants $c_i\neq c_j$, then $\cF_i\cup \cF_j$ is not a global minimum for FP (in  the adequate restriction), since by the structure given in Thm. \ref{T:minimos globales}, if a global minimum is a union of two mutually orthogonal tight subframes (on their spans),  then one of them must be a single vector. 
On the other side, if every  $\cF_i$, consisting of more than one vector has the same frame constant $c$, then there must be a $j$ such that $\cF_j$ has only a single vector, with $c_j\neq c$ (since $\cF$ can not be a tight frame).  Moreover, by Remark \ref{ordenadito} and Thm. \ref{T:minimos globales}, $c_j<c$ which implies that $\cF_j\cup \cF_i$ is not a global minimum, again by Thm. \ref{T:minimos globales}. 

So, let $\cF_i$, $\cF_j$ be such pair of subsets. Notice that if $c_i>c_j$, then the vectors in $\cF_j$ must be linear dependent, since it always have more than one vector (recall that the partition on orthogonal subsets of $\cF$ is minimal). Then, from the proof of Claim 3 in the proof of \cite[Thm. 10]{casazza2} we deduce that given $\eps>0$ there exist a set $\cF(\eps)$ such that $d(\cF_i\cup\cF_j,\cF(\eps))\leq \eps$ and FP$(\cF(\eps))<$FP$(\cF)$.	 
Hence, $\cF$ is not a local minimizer of FP on $\cB(\a)$.
\end{proof}

By inspection of the previous proof, we see that the complete characterization of local minimum for every $P_f$ on $\cB(\a)$ depends on the following problem:

\medskip

{\bf Problem($\star$):} let $\cF=\cF_1\cup \cF_2\in \cB(\a)$ such that $\cF_1\perp \cF_2$ and $\cF_i$ is a tight frame on its span. Suppose that $\cF$ is not a global minimum for $P_f$. Given $\eps>0$. Is there a frame $\cG\in \cB(\a)$ such that $\|\psi_i-\phi_i\|<\eps, \forall i$ and $P_f(\cG)<P_f(\cF)$?

\section{APPENDIX: A Geometrical approach to the frame perturbation problem}

We now consider some well known facts from differential geometry that we shall need in the sequel. In what follows we consider the unitary group $\U(m)$ together with its natural differential geometric (Lie) structure. Given $U\in \U(m)$ we shall identify its tangent space $$\mathcal T_U\,\U(m)=\{X\in \mathcal M_m(\CC):\ U^*X\in i\cdot\mathcal M_m(\CC)^{sa}\}$$ with the fixed space $\mathcal T_I\,\U(m)=i\cdot \mathcal M_m(\CC)^{sa}$ of $m\times m$ anti-hermitian matrices, via the isometric isomorphism $X\mapsto U^*X$. Given $G\in \mathcal M_m(\CC)^+$ we consider the smooth map $\Psi_G:\U(m)\rightarrow \U_m(G)$ given by $\Psi_G(U)=U^*GU$. Under the previous identification of the tangent spaces of $\U(m)$, the differential of $\Psi_G$ at a point $U\in \U(m)$ in the direction given by $X\in i\cdot \mathcal M_m(\CC)^{sa}$ is given by \begin{equation}\label{prediffi}
(D\Psi_G)_U(X)=[X,U^*GU].
\end{equation} As it is well known, the differential $(D\Psi_G)_U$ is an epimorphism at every $U\in \U(m)$ and hence \eqref{prediffi} gives us a description of the tangent space of the manifold $\U_m(G)$ at a point $U^*GU$.

Let $\Delta(G)=\{x\in \RR^m:\ \sum_{i=1}^m x_i=\tr(G) \}$ and consider $\Phi_G:\U(m)\rightarrow \Delta(G)$ given by $\Phi_G(U)=\text d(U^*GU)$, where $\text d(A)\in \RR^m$ is the main diagonal of the matrix $A\in \mathcal M_m(\CC)$. Notice that $\Delta(G)$ is a sub-manifold of $\RR^m$ with tangent space at $x\in \Delta(G)$ $$\mathcal T_x \Delta(G)=\{ y\in \RR^m: \ \sum_{i=1}^m y_i=0\}.$$
Using \eqref{prediffi}, we get (identifying again the tangent spaces of $\U(m)$ as before) that the differential of $\Phi_G$ at a point $U\in \U(m)$ in the direction of $X\in i\cdot \mathcal M_m(\CC)^{sa}$ is
\begin{equation}\label{prelims diffi2}
(D\Phi_G)_U(X)=\text d([X,U^*GU]).
\end{equation}
We shall be concerned with the existence of local cross sections of the map $\Phi_G$ around the identity $I\in \U(m)$. Since the map $\Phi_G$ is smooth, the existence of local cross sections of $\Phi_G$ is equivalent to the surjectivity of its differential $(D\Phi_G)_I$ around the identity. 

Let us fix some notation first:  we shall denote by $\mathbb{I}_m$ the (ordered) set $(1,2,\ldots,m)$. Let $\{e_i\}_{i\in \mathbb{I}_m}$ be the canonical orthonormal basis in $\CC^m$, for $I\subseteq \mathbb{I}_m$ we let $P_I$ denote the (diagonal) projection onto the span$\{e_i:\ i\in I\}$. Finally,  by $B_\delta(x)$ we mean a ball centered on $x$ with radius $\delta$, in the metric given by the context. 

The following result is part of Step 1 in \cite{leite}.

\medskip

\medskip

\begin{lem}\label{cond nec para sec loc}
Let $G\in \mathcal M_m(\CC)^+$ with $\text d(G)=\a$ and consider $\Phi_G$ as before. Then, the differential $(D\Phi_G)_I:i\cdot \mathcal M_m(\CC)^{sa}\rightarrow \mathcal T_\a\Delta(G)$ is surjective, and hence $\Phi_G$ is open in $\Delta(G)$,  if for $I\subseteq \mathbb I_m$ such that $P_IG=GP_I$ then $I=\mathbb I_m$ or $I=\emptyset$.
\end{lem}
\begin{proof}
Assume that $(D\Phi_G)_I$ is not surjective. Then, there exists $0\neq x\in \mathcal T_\a\Delta(G)$ which is orthogonal to the image of $(D\Phi_G)_I$. Let $D$ be the diagonal matrix with main diagonal $x\in \RR^m$. Using \eqref{prelims diffi2} we get  
\begin{equation}\label{equadif}
0=\langle \text d([X,G]),\, x\rangle=\tr( [X,G]D)=\tr(X [G,D]), \quad \forall X\in i\cdot \cM_m(\CC)^{sa}.
\end{equation}Since $[G,D]$ is also anti-hermitian we get that $[G,D]=0$ and hence $G$ and $D$ commute. If we let $I=\{ i: x_i>0\}$ we see, since $P_I$ is a polynomial in $D$, that $[G,P_I]=0$. Notice that $I\neq \emptyset$ and $I\neq \mathbb I_m$ since $\sum_{i=1}^m x_i=0$.
\end{proof}

\begin{lem}\label{jugadetti}
 Let us assume that the map $\Phi:=\Phi_G$, defined as before for $G\in \mathcal M_m(\CC)^+$, has a local cross section around the identity. Let $\{ G_i\}_i\subseteq  \mathcal M_m(\CC)^+$ be a sequence converging to $G$ and for $i\in \NN$ let $\Phi_i:=\Phi_{G_i}$ be defined as before. 
Then there exist $\delta>0$ and $i_0\in \NN$ such that for $i\geq i_0$ then
 $$B_\delta(I)\cap \U(m)=\mathcal S+\mathcal K_i$$ where $\mathcal S$ and $\mathcal K_i$ are submanifolds with $ I=(I_\mathcal S,I_{\mathcal K_i})$ and $$\Phi_i|_{\mathcal S}:\mathcal S\rightarrow \Phi_i(\mathcal S)\, ,\ \ \ \Phi|_{\mathcal S}:\mathcal S\rightarrow \Phi(\mathcal S)$$ are diffeomorphisms.
\end{lem}
\begin{proof}
First note that without loss of generality we can assume, as we shall, that $\tr(G_i)=\tr(G)$ for $i\in \NN$. Also note that the maps $\Phi_i$ converges uniformly to $\Phi$ since  
\begin{equation}\label{conv unif} 
\Phi_i(U)-\Phi(U)=\text d(U^*(G_i-G)U).
\end{equation} 

On the other hand, there is uniform convergence at the level of the differentials of these transformations. Indeed, under the previous identification of the tangent spaces of $\U(m)$ we can apply \eqref{prelims diffi2} and get 
\begin{equation}\label{dif unif}
\|(D\Phi)_U(X)-(D\Phi_i)_U(X)\|=\|\text d([X,U^*(G-G_i)U])\|\leq 2\sqrt{m} \,\|X\| \,\|G-G_i\|.
\end{equation} where $X\in i\cdot \mathcal M_m(\CC)^{sa}$ is arbitrary. 

We now consider $\Gamma:W\rightarrow B_{\delta_1}(I)\cap \U(m)$ a diffeomorphic local chart, where $W\subseteq \RR^p$ is an open set with $\Gamma(0)=I$. Let $\Phi\circ \Gamma:W\rightarrow \Delta(G)$ and notice that $(D (\Phi \circ \Gamma))_0:\RR^p\rightarrow \mathcal T_\a\Delta(G
)$ is surjective.  By continuity, we can assume that $(D(\Phi\circ \Gamma))_x$ is surjective for all  $x\in W$. Hence, the orthogonal projection $Q_x$ to $(\ker (D (\Phi\circ \Gamma))_x)^\perp$ is continuous on $W$. Indeed in this case we have that $Q_x=D_x^*(D_x D_x^*)^{-1}D_x$ since $D_x:=(D(\Phi\circ \Gamma))_x$ is surjective on $W$. By continuity of the projections $Q_x$ we can assume without loss of generality that $\|Q_0(1-Q_x)\|\leq 1/4$ for all $x\in W$.

By taking $0<\delta\leq\delta_1$ and using the uniform convergence of the differentials \eqref{dif unif}, we can assure that there exists a $i_1\in\NN$ such that for all $i\geq i_1$ then $(D(\Phi_i\circ\Gamma))_x$ is surjective for all $x\in W$. If $Q_{x,\,i}$ denotes the orthogonal projection onto $(\ker (D(\Phi_i\circ \Gamma))_x)^\perp$ then, using the previous description of $Q_{x,\,i}$ we see that for every $\epsilon>0$ there exists $i(\epsilon)$ such that $\|Q_{x,\,i}-Q_x\|\leq \epsilon$ for $i\geq i(\epsilon)$ and for every $x\in W$. Let $i_2=i(1/4)\in \NN$, then if $i_0=\max\{i_1,1_2\}$,  for every $x\in W$ and every $i\geq i_0$ we have
$$
\|Q_0(1-Q_{x,\,i})\|\leq \|Q_0(Q_x-Q_{x,\,i}) \|+\|Q_0(1-Q_x)\|\leq 1/2
$$
and hence \begin{equation}\label{no se tocan} (\ker (D( \Phi \circ\Gamma ))_0)^\perp\cap\ker (D (\Phi_i\circ\Gamma ))_x=\{0\}.\end{equation}
We now define $\mathcal S:=\Gamma( \ker (D (\Phi\circ \Gamma))_0)^\perp \cap W)$ 
and $\mathcal K:=\Gamma(\ker (D(\Phi\circ \Gamma))_0 \cap W)$. An straightforward argument using \eqref{no se tocan} now shows that $(D \Phi|_\mathcal S)_x$ is injective and using a dimension argument we conclude that $(D \Phi|_\mathcal S)_x$ is also surjective for all $x\in \mathcal S$; similarly with $\Phi_i$ for $i\geq i_0$. The lemma follows from these last facts.

\end{proof}

\begin{lem}
Using the notations and assumptions of the previous lemma, let $\Psi:A\,(=A^0\subseteq \RR^t)\rightarrow \mathcal S$ be a local chart of $\mathcal S$ with $\Psi(0)=I_\mathcal S$ and let $V(r):=\Psi(\text B_r(0))\subseteq \mathcal S$, where $\overline {\text B_r(0)}\subseteq A$.  Then, for any such $r>0$ there exists $\varepsilon>0$ such that  for $i\geq i_0$ then
\begin{equation}\label{albiertito}
\text B_\varepsilon(\text d(G_i))\subseteq \Phi_i(V(r)).
\end{equation}
\end{lem}
\begin{proof}
Fix $r$ as above and let $V=V(r)$. Note that for $i\geq i_0$ then $\Phi_i(I_\mathcal S)=\text d(G_i)$ is an interior point of $\Phi_i(V)$ and similarly $\Phi(I_\mathcal S)=\text d(G)$ is an interior point of $\Phi(V)$. 
We show that there exists $\epsilon >0$ such that for all $i\geq i_0$ 
then \begin{equation}\label{dist}
\inf_{x\in \partial \Phi_i(V) }\| \text d(G_i)-x\|=\min_{x\in \partial \Phi_i(V) }\| \text d(G_i)-x\|\geq \epsilon
\end{equation} where $ \partial \Phi_i(V) $ stands for boundary of the image $\Phi_i(V)$ in $\Delta(G)$. Observe that the lemma is a consequence of the condition given in \eqref{dist}.

Indeed, assume that \eqref{dist} is not true. Then, there exists a (sub)-sequence $(\Phi_{i_k})$ such that 
\begin{equation}\label{todo mal}
\inf_{x\in \partial \Phi_{i_k}(V) }\| \text d(G_{i_k})-x\|=\|\ \text d(G_{i_k})-x_k\| \leq \frac{1}{k}
\end{equation} for some $x_k=\Phi_{i_k}(U_k)$ with $U_k\in \partial V\subseteq \mathcal S$ since 
$\Phi_{i_k}(\partial V)=\partial \Phi_{i_k}(V)\subseteq \Delta$. But then for every $k\in \NN$ and $U_k\in \partial V$ then 
\begin{eqnarray*}\label{la ultima por hoy}
\|\text d(G)- \Phi(U_k)\|&\leq& \|\text d(G)-\text d(G_{i_k})\|+\|\text d(G_{i_k})-\Phi_{i_k}(U_k)\|+
\| \Phi_{i_k}(U_k)- \Phi(U_k)\| \\ &=& \|\text d(G)-\text d(G_{i_k})\|+\|\text d(G_{i_k})-x_k\|+
\| \Phi_{i_k}(U_k)- \Phi(U_k)\|\xrightarrow [k]{}0
\end{eqnarray*}by \eqref{todo mal} and the convergences  $\text d(G_{i_k})\rightarrow \text d(G)$ and 
$\Phi_{i_k}(U_k)\rightarrow \Phi(U_k)$. But this implies that $\text d(G)$ is not an interior point of $\Phi(V)$ since in this case $$\inf_{x\in \partial \Phi(V)}\|\text d(G)-x\|=\inf_{z\in\partial V}\|\text d(G)-\Phi(z)\|=0$$ which contradicts the claims at the beginning of this proof.
\end{proof}

\begin{teo}\label{anda con sec loc}
Let $\cF=\{\phi_j\}_{j=1}^m\subseteq \CC^d$ be a list of vectors, let $G=G^\cF$ be its Grammian operator and assume that $\Phi:=\Phi_G$ has a local cross section around the identity. 

Let $\{S_i\}_i\in \cM_d(\CC)^+$ be a sequence converging to $S=S^\cF$. 
Then, for every $\eta>0$ there exists $i_1\in \NN$ such that for each $i\geq i_1$ there exists a frame $\cG(i)=\cG=\{ \psi_j\}_{j=1}^m$ such that: 

\begin{enumerate}
\item\label{itt1} $\|\psi_j\|=\|\phi_j\|$ for $1\leq j\leq m$.
\item\label{itt2} $\|\psi_j-\phi_j\| \leq \eta$ for $1\leq j\leq m$.
\item\label{itt3} $S^\cG=S_i$.
\end{enumerate}
\end{teo}
\begin{proof}
Let $T=T^\cF:\CC^m\rightarrow \CC^d$ be the frame operator of the list $\cF$ with polar decomposition 
$T=|T^*|\,W=S^{1/2}\, W$ for a co-isometry $W:\CC^m\rightarrow \CC^d$. Define $G_i=W^*S_iW$ and notice that, by our hypothesis, $\|G_i-G\|\xrightarrow[]{i} 0$.

Using the notation introduced in the previous lemmas, let $\Psi:A\,(=A^0\subseteq \RR^t)\rightarrow \mathcal S$ be a local chart and $r>0$ be small enough so that $B_r(0)\subseteq A$ and for $U\in V(r)=\Psi(B_r(0))$ then \begin{equation}\label{acot1}\|U-I\|\leq \frac{\eta}{2(\|S^{1/2}\|+\eta/2)}.\end{equation}
For this choice of $r>0$ let $\varepsilon>0$ be as in \eqref{albiertito} for $i\geq i_0\in\NN$.
Let $i_2\in \NN$ be such that, for $i\geq i_2$ then  $\|S^{1/2}-S_i^{1/2}\|\leq \eta/2$ and $\|G_i-G\|\leq \frac{\varepsilon}{\sqrt{m}}$. 

If we now define $i_1=\max(i_0,i_2)$ then for $i\geq i_1$ we further have 
\begin{equation}\label{acot 3} \|\Phi(I)-\Phi_i(I)\|=\|\text d(G-G_i)\|\leq \sqrt{m}\,\|G-G_i\|< \varepsilon \ \Rightarrow \ \text d(G)\in \Phi_i(V(r)).
\end{equation}

We fix $i\geq i_1$ and construct $\cG=\cG(i)$ with the desired properties. By Lemma \ref{jugadetti} and \eqref{acot 3} there exists $U\in V(r)\subseteq \mathcal S$ such that $\Phi_i(U)=\text d(G)$. 

Define $\tilde T :=S_i^{1/2}\, W U$, and $\cG =\{\psi_j\}_{j=1}^m=\{\tilde T(e_j)\}_{j=1}^m$ where $\{e_j\}_{j=1}^m$ denotes the canonical basis of $\CC^m$. Since by construction $G^\cG=U^*G_iU$  and $S^\cG=S_i$, then items  \eqref{itt1} and \eqref{itt3} hold true. Item \eqref{itt2} follows from the inequality 
$$\|T-\tilde T\|=\|S^{1/2}W-S_i^{1/2}WU\|\leq \|S^{1/2}-S_i^{1/2}\|+\|S_i^{1/2}\|\, \|I-U\|\leq \eta.$$

\end{proof}

\begin{proof}[ Proof  of Proposition \ref{prop: anda con sec loc}] This is an immediate consequence of Lemma \ref{cond nec para sec loc} and Theorem \ref{anda con sec loc}.
\end{proof}

\end{document}